\documentclass[aos,MSNbibl,nameyear,seceqn,dvips]{arximspdf}

\doi{10.1214/13-AOS1108} 
\volume{41}
\issue{3}
\pubyear{2013}
\firstpage{1260}
\lastpage{1267}

\makeatletter
\newcommand{\eqref}[1]{(\ref{#1})}
\def\cal{\mathcal}

\newcommand{\R}{\mathbb{R}}
\newtheorem{satz}{Theorem}[section]
\newtheorem{lem}[satz]{Lemma}
\newproclaim{defin}[satz]{Definition}

\newproclaim{exam}[satz]{Example}
\makeatother

\begin{document}
\begin{frontmatter}

\title{Complete classes of designs for nonlinear regression models and
principal representations of~moment spaces\thanksref{T1}}
\thankstext{T1}{Supported in part by
the Collaborative Research Center ``Statistical modeling of
nonlinear dynamic processes'' (SFB 823, Teilprojekt C2) of the
German Research Foundation (DFG).}
\runtitle{Complete classes}

\begin{aug}
\author{\fnms{Holger} \snm{Dette}\corref{}\ead[label=e1]{holger.dette@rub.de}}
\and
\author{\fnms{Kirsten} \snm{Schorning}}

\runauthor{H. Dette and K. Schorning}

\affiliation{Ruhr-Universit\"at Bochum}
\address{Fakult\"at f\"ur Mathematik\\
Ruhr-Universit\"at Bochum\\
44780 Bochum\\
Germany\\
\printead{e1}}
\dedicated{In memory of W. J. Studden}
\end{aug}

\received{\smonth{11} \syear{2012}}
\revised{\smonth{2} \syear{2013}}

%
\begin{abstract}
In a recent paper Yang and Stufken [\textit{Ann. Statist.} \textbf{40}
(2012a) 1665--1685] gave sufficient
conditions for complete classes of designs for nonlinear regression
models. In this note we
demonstrate that there is an alternative way to validate this result.
Our main argument utilizes
the fact that boundary points of moment spaces generated by Chebyshev
systems possess unique representations.
\end{abstract}




\end{frontmatter}

\section{Introduction}
\label{sec1}

The construction of locally optimal designs for nonlinear regression
models has found considerable interest in recent years [see, e.g., \citet
{hestusun1996}, \citet{detmelwon2006}, \citet{khumuksingho2006}, \citet
{fangheda2008}, \citet{yangstuf2012a} among others]. While most of the
literature focuses on specific models or specific optimality criteria,
general results characterizing the structure of locally optimal designs
are extremely difficult to obtain due to the complicated structure of
the corresponding nonlinear optimization problems. In a series of
remarkable papers \citet{yangstuf2009}, \citet{yang2010}, \citet
{detmel2011} and \citet{yangstuf2012b} derived several complete classes
of designs with respect to the Loewner Ordering of the information
matrices. The first paper in this direction of \citet{yangstuf2009}
investigates nonlinear regression models with two parameters. These
results were generalized by \citet{yang2010} and \citet{detmel2011} to
identify small complete classes for nonlinear regression models with
more than two parameters. The most general contribution is the recent
paper of \citet{yangstuf2012b}, which provides a sufficient condition
for a complete class of designs and is applicable to most of the
commonly used regression models. On the one hand, the proof of this
statement is self-contained and only involves basic algebra. On the
other hand, the proof is complicated, requires several auxiliary
results and hides
some of the mathematical structure of the problem.

The purpose of the present paper is to demonstrate that conditions of
this type are intimately related to the characterization of boundary
points of moment spaces associated with a nonlinear regression model.
Our main tool is a Chebyshev system [\citet{karstu1966}] appearing in (a~transformation of) the Fisher information matrix of a given design. The
complete class of designs can essentially be characterized as the set
of measures corresponding to the unique representations of the boundary
points of the corresponding moment spaces. With this insight the main
result in the paper of \citet{yangstuf2012b} is a simple consequence of
the fact that a representation of a boundary point of a
$k+1$-dimensional moment space associated with a Chebyshev system
depends only on the first $k$ functions which are used to generate the
moment space.

In Section \ref{sec2} we state some facts about
moment spaces associated with Chebyshev systems which are of general
interest for constructing admissible designs. The design problem and
Theorem 1 of \citet{yangstuf2012b} are stated in Section~\ref{sec3},
where we also present our alternative proof. We finally note that the
paper of \citet{yangstuf2012b} contains numerous interesting examples
and provides a further result which are not discussed in this note for
the sake of brevity.

\section{Chebyshev systems and associated moment spaces}
\label{sec2}
A set of $k$ real valued functions $\Psi_0,\ldots,\Psi_{k-1}\dvtx [A,B]
\rightarrow\mathbb{R}$ is called {Chebychev system} on the interval
$[A,B]$ if and only if it fulfills the inequality
\[
\det \pmatrix{ \Psi_0(x_0)&\ldots&
\Psi_0(x_{k-1})\vspace*{2pt}
\cr
\vdots& \ddots& \vdots
\vspace*{2pt}
\cr
\Psi_{k-1}(x_0)&\ldots&
\Psi_{k-1}(x_{k-1})} >0
\]
for any points $x_0,\ldots,x_{k-1}$ with $A \leq x_0<x_1<\cdots<x_{k-1}
\leq B$.
The moment space associated with a Chebyshev system is defined by
\begin{eqnarray*}
\mathcal{M}_{k-1}&=& \biggl\lbrace c=(c_0,
\ldots,c_{k-1})^T \Big| c^0_i=\int
_A^B\Psi_i(x)\,d\sigma(x),
\\
&&\hspace*{98pt} i=0,\ldots,k-1, \sigma\in\mathbb{P}\bigl([A,B]\bigr) \biggr\rbrace,
\end{eqnarray*}
where $\mathbb{P}([A,B])$ denotes the set of all finite measures on the
interval $[A,B]$. It can be characterized as the smallest convex cone
containing the curve
\[
\mathcal{C}_{k-1}= \bigl\lbrace\bigl(\Psi_0(t), \ldots,
\Psi_{k-1}(t)\bigr)^T | t\in[A,B] \bigr\rbrace;
\]
see \citet{karstu1966}.
By Caratheodory's theorem, any point of $\mathcal{M}_{k-1}$ can be
described as a linear combination of at most $k+1$ points in $\mathcal
{C}_{k-1}$, where the coefficients are positive. Moment spaces can be
defined for any set of linearly independent functions, but if the
functions $\{\Psi_0, \ldots, \Psi_{k-1}\}$ generate a Chebyshev system,
the moment space has several additional interesting properties. In
particular, fewer points of $\mathcal{C}_{k-1}$ are required for the
representation of points in $\mathcal{M}_{k-1}$.
To be precise, we define for a point $c^0\in\mathcal{M}_{k-1}$ its
index $I(c^0)$ as the minimal number of points in $\mathcal{C}_{k-1}$
which are required to represent $c^0$, where the points $(\Psi_0(A),
\ldots, \Psi_{k-1}(A))^T$ and $(\Psi_0(B), \ldots, \Psi_{k-1}(B))^T$
corresponding to the boundary points of the interval $[A,B]$ are
counted by $1/2$. The index $I(\sigma)$ of a finite measure $\sigma$ on
the interval $[A,B]$ is defined as the number of its support points,
where the boundary points are counted as $1/2$. If $c^0=\int_A^B(\Psi
_0(x), \ldots, \Psi_{k-1}(x))^T\,d\sigma(x)$, the measure $\sigma$ is
also called a representation of the point $c^0 \in\mathcal{M}_{k-1}$.
If $\{ t_1,\ldots,t_n \}$ denotes the support of $\sigma$, the vectors
$\{ ( \Psi_0(t_j), \ldots, \Psi_{k-1}(t_j))^T \mid j=1,\ldots,n \}$ and
the corresponding weights of $\sigma$ can be used to obtain a convex
representation of the $c^0$ by elements of $\mathcal{C}_{k-1}$.

With this convention it follows that the point $c^0\in\mathcal
{M}_{k-1}$ is a boundary point of $\mathcal{M}_{k-1}$ if and only if
its index satisfies $I(c^0) < \frac{k}{2}$.
Similarly, $c^0$ is in the interior of $\mathcal{M}_{k-1}$ if its index
is $\frac{k}{2}$.
Following \citet{karstu1966} we denote a representation $\sigma$ of an
interior point $c^0$ as \textit{principal}, if $I(\sigma)=I(c^0)=\frac{k}{2}$.
These authors also proved that representations of boundary points are
unique. Furthermore, for each interior point $c^0 \in\mathcal
{M}_{k-1}$ there exist exactly two principal representations (a~further
proof of this statement is given below). The first is called \emph
{upper} principal representation and contains the point $B$ of the
interval $[A,B]$, whereas the second is called \emph{lower} principal
representation and does not use this point. These measures are denoted
by $\sigma^+$ and $\sigma^-$, respectively. If $k$ is odd, the lower
and upper principal representation has $\frac{k+1}{2}$ support points.
On the other hand, if $k$ is even, the lower and upper principal
representation have $\frac{k}{2}$ and $\frac{k+2}{2}$ support points,
respectively. The next Lemma is crucial in the following investigations.

\begin{lem} \label{lem1}
Let $\Psi_j\dvtx [A,B] \to\mathbb{R}$ $(j=0,\ldots,k-1);  \Omega\dvtx  [A,B]
\to\mathbb{R}$ denote real valued functions and assume
that the systems $\{ \Psi_0,\ldots,\Psi_{k-1}\}$ and $\{\Psi_0,\ldots, \Psi_{k-1},\Omega\}$ are {Chebyshev systems} on the interval $[A,B]$.
If $c^0=(c^0_1, \ldots, c^0_{k-1})^T \in{\cal M}_{k-1}$,
then the upper and lower principal representation $\sigma^+ $ and
$\sigma^-$ of $c^0$ are uniquely determined and satisfy
\begin{eqnarray*}
&&\max \biggl\lbrace\int_A^B\Omega(t)\,d\sigma(t)
\bigg| \sigma\in\mathbb {P}\bigl([A,B]\bigr), c^0_i=\int
_A^B\Psi_i(t)\,d\sigma(t), i =0,
\ldots, k-1 \biggr\rbrace
\\
&&\qquad= \int_A^B\Omega(t)\,d\sigma^+(t),
\\
&&\min \biggl\lbrace\int_A^B\Omega(t)\,d\sigma(t)
\bigg| \sigma\in\mathbb {P}\bigl([A,B]\bigr), c^0_i=\int
_A^B\Psi_i(t)\,d\sigma(t), i =0,
\ldots, k-1 \biggr\rbrace
\\
&&\qquad= \int_A^B\Omega(t)\,d\sigma^-(t).
\end{eqnarray*}
In particular both representations do not depend on the function \mbox{$\Omega
\dvtx [A,B] \rightarrow\mathbb{R}$}.
\end{lem}

\begin{pf}The proof follows essentially from the discussion in
Sections 3--5 of Chapter II in \citet{karstu1966} and---as proposed by
a referee---some details are given here for sake of completeness.
If $c^0$ is a boundary point of the moment space
$\mathcal{M}_{k-1}$, there exists precisely one representation, say
$\sigma^0$, of $c^0$. This shows that the
set of measures $\sigma\in\mathbb{P} ([A,B])$ satisfying $c^0_i = \int^B_A \Psi_i(x) \,d \sigma(x)\ (i=0,\ldots,k-1)$ is a singleton, which
yields $\sigma^0= \sigma^+ = \sigma^-$ and the statement of Lemma \ref
{lem1} is obvious.

Therefore it remains to consider the case where $c^0$ is an interior
point of the moment space $\mathcal{M}_{k-1}$, that is, $I(c^0)=\frac
{k}{2}$. We assume that $k=2m$ and that there exist two upper principal
representations, say $\sigma^+_1$ and $\sigma^+_2$ (the case $k=2m-1$
and the corresponding statement for the lower principal representation
are shown by similar arguments). Because $I(\sigma^+_1)=I(\sigma
^+_2)=I(c_0)=m$, it follows that $\sigma^+_1$ and $\sigma^+_2$ have
$m+1$ support points including the boundary points $A$ and $B$. Now, if
$\sigma^+_1 \neq\sigma^+_2$, the signed measure $\sigma^+_1 - \sigma
^-_1$ has at most $2m$ support points and satisfies
\[
0= \int^B_A \bigl(\Psi_0(x),
\ldots,\Psi_{2m-1}(x)\bigr)^T \,d\bigl(\sigma^+_1 -
\sigma ^+_2\bigr) (x).
\]
Because $\{ \Psi_0, \ldots, \Psi_{2m-1} \}$ is a Chebyshev system, it
follows that $\sigma^+_1 = \sigma^+_2
$, which proves the first part of Lemma \ref{lem1}.

For a proof of the second part we note that the set
\[
\biggl\lbrace\int_A^B\Omega(t)\,d\sigma(t) \bigg|
\sigma\in\mathbb {P}\bigl([A,B]\bigr), c^0_i=\int
_A^B\Psi_i(t)\,d\sigma(t), i =0,
\ldots, k-1 \biggr\rbrace
\]
is a bounded closed interval, say $[\gamma^-, \gamma^+]$. Moreover, the
points $c^-_0=(c^T_0, \gamma^-)^T$ and $c^+_0=(c^T_0, \gamma^+)^T$ are
boundary points of the moment space $\mathcal{M}_{2m}$ generated by the
Chebyshev system
\[
\{ \Psi_0, \ldots, \Psi_{2m-1}, \Omega\}.
\]
Consequently, $I(c^{\pm}_0) < \frac{2m+1}{2}$ and the representations
of $c^+_0$ and $c^-_0$ are unique. Moreover, because $I(c_0)=m$ we also
have $I(c^{\pm}_0)=m$. {It is shown in Karlin and Studden
[(\citeyear{karstu1966}), pages
55--56]} that the representations of $c^+_0$ and $c^-_0$ must coincide
with the principal representations $\sigma^+$ and $\sigma^-$ of the
interior point $c^0 \in\mathcal{M}_{k-1}$, which proves the second
assertion of Lemma \ref{lem1}.
\end{pf}

\section{A complete class of designs for regression models}
\label{sec3}

Consider the common nonlinear regression model
%
\begin{equation}
E[Y|x] =\eta(x, \theta),
\end{equation}
where $\theta\in\mathbb{R}^p$ is the vector of unknown parameters,
$x$ denotes a real valued covariate from the design space $[A,B]\subset
\mathbb{R}$ and different observations
are assumed to be independent with
variance $\sigma^2$. The function $\eta$ is called regression
function [see \citet{sebwil1989} or
\citet{ratkowsky1990}] and assumed to be continuous and differentiable
with respect to
the variable $\theta$. A design is defined as a probability measure $\xi
$ on the interval
$[A,B] $ with finite support; see \citet{kiefer1974}. If the design $\xi
$ has masses $w_i$ at the points $x_i\ (i = 1, \ldots, l)$ and $n$
observations can be made by the experimenter, this means that the
quantities $w_i n$ are rounded to integers, say $n_i$, satisfying
$\sum^l_{i=1} n_i =n$, and the experimenter takes $n_i$ observations at
each location $x_i\ (i=1, \ldots, l)$.
If the design $\xi$ contains $l$ support points $x_1,\ldots, x_l$ such
that the vectors $\frac{\partial}{\partial\theta}\eta(x_1,\theta),
\ldots, \frac{\partial}{\partial\theta}\eta(x_l,\theta)$ are linearly
independent, and observations are taken according to
this procedure, it follows from \citet{jennrich1969} that the covariance
matrix of the
nonlinear least squares estimator is approximately (if $n \to\infty$)
given by
%
\begin{equation}
\label{var} {\sigma^2\over n} M^{-1}(\xi, \theta)=
{\sigma^2\over n} \biggl(\int_A^B \biggl(
\frac
{\partial}{\partial\theta}\eta(x,\theta) \biggr) \biggl(\frac{\partial
}{\partial\theta}\eta(x,\theta)
\biggr)^T\,d\xi(x) \biggr)^{-1}.
\end{equation}
An optimal design maximizes an appropriate functional of the
matrix\break
${n\over\sigma^2}
M(\xi, \theta)$, and numerous criteria have been proposed
in the literature to discriminate between competing designs; see \citet
{pukelsheim2006}. Note that the matrix \eqref{var}
depends on the unknown parameter $\theta$, and following \citet
{chernoff1953} we call the maximizing designs
locally optimal designs. These designs require an initial guess of the
unknown parameters in the model and are used as benchmarks for many
commonly used designs or for the construction of more sophisticated
optimality criteria which require less information regarding the
parameters of the model [\citet{chaver1995} and \citet{dette1997}].

Most of the available optimality criteria are positively homogeneous, that
is, $
\Phi ({n\over\sigma^2}M(\xi, \theta) ) = {n\over\sigma^2} \Phi
(M(\xi, \theta))
$ [\citet{pukelsheim2006}]. Therefore it is sufficient to consider
maximization of functions of the matrix $M(\xi, \theta)$, which is
called \textit{information matrix} in the literature.
Moreover, the commonly used optimality criteria also satisfy a
monotonicity property with respect to the Loewner ordering, that is,
$ \Phi(M(\xi_1, \theta)) \geq\Phi( M(\xi_2, \theta)),
$
whenever $M(\xi_1, \theta) \geq M(\xi_2, \theta)$,
where the parameter $\theta$ is fixed, $\xi_1$, $\xi_2$ are two
competing designs on the interval $ [A,B]$
and $\Phi$ denotes an information function in the sense of \citet
{pukelsheim2006}.
Throughout this paper we call a design $\xi$ admissible if there does
not exist any design $\xi_1$, such
that $M(\xi_1, \theta)\neq M(\xi, \theta)$ and
%
\begin{equation}
\label{infoleq} M(\xi_1, \theta) \geq M(\xi, \theta).
\end{equation}
\citet{yangstuf2012b} derive a complete class theorem in this general
context which characterizes the class of designs, which cannot be
improved with respect to the Loewner ordering of their information
matrices. For the sake of completeness and because of its importance we
will state this result
here again.
In particular, we demonstrate that the complete class specified by
these authors corresponds to upper and lower principal representations
of a moment space generated by the regression functions.
For this purpose we denote by $P(\theta)$ a regular $p\times p$ matrix,
which does not depend on the design $\xi$, such that
the representation
%
\begin{equation}
M(\xi, \theta)= P(\theta)C(\xi, \theta)P^T(\theta)
\end{equation}
holds, where the $p\times p$ matrix $C(\xi, \theta)$ is defined by
\begin{eqnarray*}
C(\xi,\theta)&=&\int_A^B \pmatrix{
\Psi_{11}(x)& \ldots& \Psi_{1p} (x)\vspace*{2pt}
\cr
\vdots&
\ddots& \vdots\vspace*{2pt}
\cr
\Psi_{p1}(x) & \ldots&
\Psi_{pp}(x) } \,d\xi(x)
\\
&=& \int_A^B\pmatrix{ C_{11}(x) &
C_{21}^T(x) \vspace*{2pt}
\cr
C_{21}(x) &
C_{22}(x)} \,d\xi(x),
\end{eqnarray*}
and $C_{11}(x) \in\R^{p-p_1\times p-p_1}$, $C_{21}(x) \in\R^{
p_1\times p-p_1} $, $C_{22}(x) \in\R^{ p_1\times p_1}$ are appropriate
block matrices
$(1\leq p_1 \leq p)$.
Obviously, $P(\theta) $ could be chosen as identity matrix, but in
concrete applications other choices
might be advantageous; see Yang and Stufken
[(\citeyear{yangstuf2012a}), Section 4] for
numerous interesting examples. A similar comment applies to the choice
of $p_1$ which is used to represent the
matrix $C$ in a $2 \times2 $ block matrix. Note that
the inequality (\ref{infoleq}) is satisfied if and only if the inequality
%
\begin{equation}
C(\xi_1, \theta) \geq C(\xi, \theta)
\end{equation}
holds. Following \citet{yangstuf2012b} we define $\Psi_0(x)=1$,
denote the different elements among $\{\Psi_{ij}|1\leq i \leq p, j\leq
p-p_1\}$
in the matrices $C_{11}(x) $ and $C_{21}(x)$ which are not constant by
$\Psi_1, \ldots, \Psi_{k-1}$
and define for any vector $Q \in\mathbb{R}^{p_1}\setminus\{0\} $ the function
%
\begin{equation}
\label{psiQ} \Psi^Q_k(x)=Q^TC_{22}(x)Q.
\end{equation}
We are now in a position to state and prove the main result of this paper.
%
\begin{satz}[{[\citet{yangstuf2012b}]}]\label{yangsatz}
\begin{enumerate}
\item[(1)] If $\{\Psi_0,\ldots,\Psi_{k-1}\}$ and $\{\Psi_0,\ldots,\Psi
_{k-1}, \Psi^Q_k\}$ are Chebyshev systems for every nonzero vector
$Q$, then for any design $\xi$ there exists a design $\xi^{+}$ with at
most $\frac{k+2}{2}$ support points, such that $M(\xi^{+},\theta) \geq
M(\xi, \theta)$.

If the index of $\xi$ satisfies $I(\xi)<\frac{k}{2}$, then the design
$\xi^{+}$ is uniquely determined in the set
%
\begin{equation}
\label{uni} \biggl\lbrace\eta \Big| \int_A^B
\Psi_i(x)\,d\eta(x)=\int_A^B
\Psi_i(x)\,d\xi (x), i=1,\ldots,k-1 \biggr\rbrace
\end{equation}
and coincides with the design $\xi$.

If the index of $\xi$ satisfies $I(\xi)\geq\frac{k}{2}$, then the
following cases are discriminated:
\begin{enumerate}[(a)]
\item[(a)] If $k$ is odd, then the design $\xi^+$ has at most $\frac
{k+1}{2}$ support points and it can be chosen such that $B$ is a
support point of the design $\xi^+$.
\item[(b)] If $k$ is even, then the design $\xi^+$ has at most $\frac
{k+2}{2}$ support points and it can be chosen such that $A$ and $B$ are
support points of the design $\xi^+$.
\end{enumerate}
\item[(2)] If $\{\Psi_0,\ldots,\Psi_{k-1}\}$ and $\{\Psi_0,\ldots,\Psi
_{k-1}, -\Psi^Q_k\}$ are Chebyshev systems for every nonzero vector
$Q$, then for any design $\xi$ there exists a design $\xi^{-}$ with at
most $\frac{k+1}{2}$ support points, such that $M(\xi^{-},\theta)\geq
M(\xi, \theta)$.

If the index of $\xi$ satisfies $I(\xi)<\frac{k}{2}$, then the design
$\xi^{-}$ is uniquely determined in the set of measures satisfying
\eqref{uni}
and coincides with the design $\xi$.

If the index of $\xi$ satisfies $I(\xi)\geq\frac{k}{2}$, then the
following cases are discriminated:
\begin{enumerate}[(a)]
\item[(a)] If $k$ is odd, then the design $\xi^-$ has at most $\frac
{k+1}{2}$ support points and it can be chosen such that $A$ is a
support point of the design $\xi^-$.
\item[(b)] If $k$ is even, then the design $\xi^-$ has at most $\frac{k}{2}$
support points.
\end{enumerate}
\end{enumerate}
\end{satz}

\begin{pf}
We only present the proof of the first part of the theorem; the second
part follows by similar arguments.
\citet{yangstuf2012b} showed that
a design $\xi_1 $ satisfies \eqref{infoleq} if the conditions
%
\begin{eqnarray}
\label{yangbdg} \int_A^B\Psi_i(x)\,d
\xi_1(x)&=&\int_A^B
\Psi_i(x)\,d\xi(x),\qquad i= 1,\ldots,k-1,
\nonumber
\\[-8pt]
\\[-8pt]
\nonumber
\int_A^B \Psi^Q_k(x)\,d
\xi_1(x)&\geq&\int_A^B
\Psi^Q_k(x) \,d\xi(x)
\end{eqnarray}
are satisfied for all vectors $Q\neq0$. Consequently an improvement of
the design $\xi$ is obtained by maximizing the ``$k$th moment''
$\int_A^B \Psi^Q_k(x)\,d\xi_1(x)$ in the set of all designs satisfying
\eqref{yangbdg}.
If $I(\xi) < \frac{k}{2}$, then this set is a singleton and the
maximizing design $\xi_Q^+$ coincides with $\xi$. Otherwise, by Lemma
\ref{lem1} the maximizing measure $\xi_Q^+$ corresponds to the upper
principal presentation of the
moment point $(\int_A^B \Psi_0(x)\,d\xi(x), \ldots, \int_A^B \Psi
_{k-1}(x)\,d\xi(x))^T $, which does not
depend on the vector~$Q$. Finally,
assertion 1(a) or 1(b) of Theorem~\ref{yangsatz} follows from the
discussion regarding the number of support points of principal
representations given at the end of Section \ref{sec2}.
\end{pf}

\section*{Acknowledgments}
The authors thank Martina
Stein, who typed parts of this manuscript with considerable
technical expertise. We are also grateful to the referees and
Associate Editor for constructive comments on an earlier version
of this paper.


%

\printaddresses

\end{document}